\documentclass{amsart}
\usepackage{graphicx}
\usepackage{amsmath,amssymb,amsbsy,amsfonts,amsthm,latexsym,
            amsopn,amstext,amsxtra,euscript,amscd,amsthm}

\newtheorem{lem}{Lemma}

\newtheorem{thm}{Theorem}
\newtheorem{theorem}[thm]{Theorem}

\newtheorem{cor}{Corollary}

\newtheorem{hyp}{Hypothesis}

\newtheorem{exam}{Example}

\newcommand{\Z}{\mathbb{Z}}

\begin{document}

\title{On the proximity of large primes}
\author{Minjia Shi}
\address{Math Dept, Anhui Hefei University, Hefei, Anhui, PR China}
\email{smjwcl.good@163.com}
\author{Florian Luca}
\address{School of Mathematics, University of the Witwatersrand, Private Bag X3, Wits 2050, South
Africa\\
Max Planck Institute for Mathematics, Vivatgasse 7, 53111 Bonn, Germany\\
Department of Mathematics, Faculty of Sciences, University of Ostrava, 30. dubna 22, 701 03
Ostrava 1, Czech Republic}
\email{florian.luca@wits.ac.za}

\author{Patrick Sol\'e}
\address{CNRS/LAGA, University of Paris 8, rue de la Libert\'e, 93 526 Saint-Denis }
\email{sole@telecom-paristech.fr}

\begin{abstract} By a sphere-packing argument, we show that there are infinitely many pairs of primes that are close to each other for some metrics on the integers. In particular, for any numeration basis $q$,
we show that there are infinitely many pairs of primes
the base $q$ expansion of which differ in at most two digits. Likewise, for any fixed integer $t,$ there are infinitely many pairs of primes, the first $t$ digits of which are the same.
In another direction, we show that, there is a constant $c$ depending on $q$ such that for infinitely many integers $m$ there are at least $c\log \log m$ primes which differ from $m$ by at most one base $q$ digit.
\end{abstract}

\maketitle
\let\thefootnote\relax
\section{Introduction}
Recently, it has been shown that there are pairs of arbitrary large primes with difference in absolute value at most 246 \cite{poly}.
In this paper, we strive to replace the absolute value in the preceding statement by some other metric.
The main thrust is to regard an integer as the vector of coefficients of its base $q$ decomposition for some integer $q>1.$ Thus we can
embed any integer of size $<q^n$ in $\Z_q^n,$ which can be equipped with a combinatorial metric, and allows us to use coding theoretic arguments. Two metrics are of importance;
first off, the Hamming metric which counts the number of different $q$-ary digits between the $q$-expansion of two integers; next, the Rosenbloom Tsfasman metric, which is related to the $q$-adic valuation of the difference of the two integers, and has appeared in the past in questions of equidistribution \cite{DS}. We consider the nonlinear code consisting of the $q$-expansion
of primes in a given integer interval. Its size can be estimated by the Prime Number Theorem, and gives an upper bound on its minimum distance. Some surprising consequences
follow. There are large primes the $q$-expansion of which differ in at most two places (Corollary 1). For any fixed integer $t\ge 1$ there are large primes coinciding on their first $t$ digits (Corollary 2). While Corollary 2 can be proved by invoking Dirichlet Theorem on primes in arithmetic progression, our proof is certainly more elementary.

Another similar, but distinct, question is to ask how many primes are close to a large integer for the Hamming metric. This question is more difficult and requires more classical techniques from analytic number theory \cite{BaSh,May}.

The article is organized as follows. The next section collects the necessary notations and definitions. Section 3 considers how close large primes can be. Section 4 studies
how close primes are to large integers. Section 5 reviews our results and points out some open questions.
\section{Definitions and Notation}
Let $q>1$ be an integer. Denote by $\Z_q$ the ring of integers modulo $q.$ We consider words of length $n$ over $\Z_q.$ We assume that $\Z_q^n$ is equipped with a metric $d(.,.).$

{\exam The Hamming metric between $x$ and $y$ in $\Z_q^n$ is defined as the number of indices $i$ where $x_i\neq y_i.$}

The ball of radius $t$ about $x$ is defined by

$$B(x,t):=\{ y \in \Z_q^n\,:\, d(y,x) \le t\}.$$

We assume that the size of $B(x,t)$ does not depend on $x,$ and we denote by $B(t)$ this common value.

{\exam For the Hamming metric it is well-known \cite{HP} that $B(t)=\sum_{i =0}^ n{n \choose i} (q-1)^i .$}

By a {\em code} of length $n$ over $\Z_q,$ we shall mean a proper subset of $\Z_q^n.$ The following result is known as the sphere packing bound.
The proof is immediate and omitted.
{\thm \label{ham}If the balls of radius $t$ centered about elements of $C$ are pairwise disjoint, then $\vert C\vert B(t) \le q^n.$}

Given a code $C,$ the largest possible $t$ such that the hypothesis of Theorem \ref{ham} holds is called the {\em error-correcting capacity} of $C$ and denoted by the letter $e.$
Denote by $d$ the minimum distance of the code. If $d$ is integer valued, then  it can be seen, using the triangle inequality, that $e\ge \lfloor \frac{d-1}{2}\rfloor.$

We also need the $q$-adic metric on integers. If $m$ is nonzero integer let $\vert m\vert_q$ be the exponent of the largest power of $q$ that divides $m.$
This number is sometimes called the $q$-adic valuation of $m.$
By convention, let $\vert 0\vert_q=0.$
The $q$-adic distance $d_q(r,s)$ between two integers is then defined as

$$ d_q(r,s)= \vert r \ominus s\vert_q+1.$$
where $\ominus$ denote the subtraction performed (without carry) on the $q$-ary expansions of $r$ and $s.$ Thus we assume that both $r,s \le q^n,$ and we write
$$ r=r_0+r_1 q+\cdots+r_{n-1}q^ {n-1},$$
as well as
$$ s=s_0+s_1 q+\cdots+s_{n-1}q^ {n-1},$$
without assuming necessarily $r_{n-1}\neq 0,$ or $s_{n-1}\neq 0.$
Thus, with that definition, we have that $d_q(0,q+q^2)=2,$ for instance.
It can be seen that this distance coincide with the RT metric on $\Z_q^n,$ up to reversal of the order of the digits.
Indeed the distance $d_{RT}(x,y)$ between two elements of $\Z_q^n$ is usually defined (cf. \cite{DS}) as the quantity

$$d_{RT}(x,y)=\max \{ i : \, x_i \neq y_i\}.  $$

If $x=(x_0, \ldots,x_{n-1})$ denote by $Rx$ the quantity $Rx=(x_{n-1},\ldots,x_0).$ With these notations, we see that

$$d_q(x,y)=d_{RT}(Rx,Ry).$$

It is easy to check that with these conventions, we have $B(t)=q^{n-t}.$ Indeed the ball of given radius centered at the origin is

$$B(0, t)=\{q^t(x_0+x_1q+\cdots+x_{n-t-1}x^{n-t-1})\},$$
where the $x_i$'s are arbitrary elements of $\Z_q.$

\section{Proximity between primes}
Let $q>1$ be an integer. For any pairs of integers, $a,b$ we define the  {\it prime code} $P(a,b)$ as the base $q$ expansions of the primes contained in $[a,b].$
If we assume that $q^{n-1}< b \le q^n,$ we can thus view this code as a non linear code of length $n$ over an alphabet of size $q.$
(Note that $0 \notin P(a,b),$ so that $P(a,b)$ cannot be a linear code). Denote by $d[a,b]$ the minimum distance of that code.
We begin with a combinatorial sphere packing argument, which is closer to the standard coding theoretic situation.
{\lem \label{sp} Assume $d$ is an integer-valued distance. If $B(t)>\frac{q^n}{\vert P(a,b)\vert},$ then $d[a,b]\le 2t.$}
\begin{proof} If $d[a,b]> 2t,$ then the hypothesis of Theorem \ref{ham} are satisfied for $C=P(a,b).$ This contradicts the hypothesis made on $B(t).$
\end{proof}

{\thm \label{comb} Assume that $d$ is integer valued, and that, for some $n_0$ and all  $n>n_0$ we have $B(t) > n\log q.$ If $a$ is dominated by $q^n$ as $n \rightarrow \infty,$
then, for large enough $n,$ there are primes $r<s$ in the range $[a,q^n ]$ that satisfy $d(r,s)\leq 2t.$}

\begin{proof}
 By the Prime Number Theorem \cite{A}, and the assumption made on $a,$ we know that $\vert P(a,q^ n)\vert$ is asymptotically equivalent to $\frac{q^ n}{\log q^ n}=\frac{q^ n}{n\log q}$ when $n\to\infty$. The result follows then  by Lemma \ref{sp}.
\end{proof}

{\bf Remark:} In the above theorem, one may take $a= \lceil \sqrt{n} \rceil,$ or even $a= \lceil \frac{n}{\log n}\rceil$ for instance. Taking $a$ a constant in $n$ would not be enough,
as we could not then exclude the case where both $r$ and $s$ are bounded above.

{\cor For any numeration basis $q$,
 there are infinitely many pairs of primes
the base $q$ expansion of which differ in at most two digits.}

\begin{proof}
We apply Theorem \ref{comb} to the case of $t=1$ and of $d$ the Hamming distance. In that case $B(1)=n(q-1)+1,$ and since $(q-1)>\log q $ for $q\ge2,$ we see that
$B(1)> n \log q$ for $n$ large enough.
\end{proof}

An alternative proof, also based on the PNT was shown to us by Igor Shparlinski.

\begin{proof}
Let $k$ be any large positive integer, let $N:=q^k$ and consider primes $p$ in the interval $[N,qN)$. The number of them is
$$
\pi(qN)-\pi(N)=(q-1+o(1))N/\log N\quad {\text{\rm as}}\quad k\to\infty,
$$
by the Prime Number Theorem. For each one of these, change one of its digits one at a time into any other digit except for the leading digit which is changed into any other digit except $0$. Since such primes have $k+1=(1+o(1))\log N/\log q$ digits as $k\to\infty$ and
$q$ digits available in each location except for the leading location where we have only $q-1$ possibilities, we get a totality of
$$
((q-1)(k+1)+O(1))(\pi(qN)-\pi(N))=\left(\frac{(q-1)}{\log q}+o(1)\right)(q-1)N
$$
as $k\to\infty$ positive integers all in the interval $(N,qN)$. Since the above interval contains $(q-1)N$ integers, and $(q-1)/\log q>1$, it follows that not all the above numbers are distinct. Thus, there are two of them which are equal. Hence, there is
an integer $m\in [N,qN)$ which arises in the above way from two distinct primes $r$ and $s$. Thus, those primes are a Hamming distance $1$ from $m$ and therefore at Hamming distance at most two from each other. So, we obtained one integer $m$ with the desired property. Since this integer was in $[q^k,q^{k+1})$, the statement about infinitely many such integers $m$ (and consequently infinitely many such pairs of primes $(r,s)$) follows by varying $k$.
\end{proof}
A concrete application of the RT metric is as follows.

{\cor For any numeration basis $q$, and any fixed integer $v$
 there are infinitely many pairs of primes $(r,s)$
such that $q^v$ divides $s- r.$}

\begin{proof}
 We apply the preceding theorem to the case when $d(.,.)=d_q(.,.)$, and $t\ge 1$ is arbitrary, and fixed. In that case $B(t)=q^{n-t}$ which
 is $> n \log q$ for $n$ large enough.
\end{proof}

{\bf Remark:} By using Dirichlet's Theorem on primes in arithmetic progressions \cite{A}, a stronger result can be proved.
Given an integer $\ell<q^v,$ and coprime with $q,$ there are infinitely many primes $ p \equiv \ell \pmod{q^v}.$
\section{Proximity to integers}
Here, we continue the previous trend and only look at positive integers $n$ which become primes by changing just one of their base $q$ digits. We prove the following theorem.

\begin{theorem}
\label{thm:1}
There exists a constant $c_q>0$ depending on $q$ such that for infinitely many $m$, $B(m,1)$ contains at least $c_q \log\log m$ prime numbers.
\end{theorem}

For this, let us recall Maynard's notion of well-distributed sets.

Given a set of integers ${\mathcal A}$, a set of primes ${\mathcal P}$, and a linear function $L(n)=l_1n+l_2$, we define
\begin{eqnarray*}
{\mathcal A}(x) &: = & \{n\in {\mathcal A}: x\le n<2x\},\\
{\mathcal A}(x;q,a) & := & \{n\in {\mathcal A}(x), n\equiv  a\pmod q\},\\
L({\mathcal A}) & := & \{L(n): n\in {\mathcal A}\},\\
 \phi_L(q) & := & \phi(|l_1|q)/\phi(|l_1|),\\
{\mathcal P}_{L,{\mathcal A}}(x) & := & L({\mathcal A}(x))\cap {\mathcal P},\\
{\mathcal P}_{L,{\mathcal A}}(x;q,a) & := & L({\mathcal A}(x;q,a))\cap {\mathcal P}.
\end{eqnarray*}
It is assumed that $({\mathcal A}, {\mathcal L}, {\mathcal P}, B,x,\theta)$ satisfy the following hypothesis where ${\mathcal L}$ is a set of linear functions of cardinality $k:=\# {\mathcal L}$, $B\in {\mathbb N}$, $x$ a large real number and $\theta\in (0,1)$:

\medskip

\begin{hyp}
\label{hip}
\begin{itemize}
\item[(i)] ${\mathcal A}$ is well distributed in progressions: we have
$$
\sum_{m\le x^{\theta}} \max_a \left|\#{\mathcal A}(x;m,a)-\frac{\#{\mathcal A}(x)}{m}\right|\ll \frac{\#{\mathcal A}(x)}{(\log x)^{100 k^2}}.
$$
\item[(ii)] Primes in $L({\mathcal A})\cap {\mathcal P}$ are well-distributed in arithmetic progressions: for any $L\in {\mathcal L}$ we have
$$
\sum_{\substack{m\le x^{\theta}\\ \gcd(m,B)=1}} \max_{\gcd(L(a),m)=1} \left|\#{\mathcal P}_{L,{\mathcal A}}(x;m,a)-\frac{\#{\mathcal P}_{L,{\mathcal A}}(x)}{\phi_L(m)}\right|\ll \frac{\#{\mathcal P}_{L,{\mathcal A}}(x)}{(\log x)^{100 k^2}}.
$$
\item[(iii)] ${\mathcal A}$ is not too concentrated on any arithmetic progression: for any $m<x^{\theta}$ we have
$$
\#{\mathcal A}(x;m,a)\ll \frac{\#{\mathcal A}(x)}{m}.
$$
\end{itemize}
\end{hyp}

Assume additionally that ${\mathcal L}=\{L_1,L_2,\ldots,L_k\}$ are such that $L_i(n)=a_i n+b_i$ satisfy $\gcd(a_i,b_i)=1$, $0\le a_i,b_i\le x^{\alpha}$ for all $i=1,\ldots,k$, $k\le (\log x)^{\alpha}$ and $B\le  x^{\alpha}$ for some fixed
$\alpha\in (0,1)$. Under these hypothesis, Maynard \cite{May} proves the following statement.

\begin{theorem}
\label{thm:Maynard}
There exists a constant $C$ depending only on $\theta$ and $\alpha$ such that if $k\ge C$, $({\mathcal A},{\mathcal L},{\mathcal P},B,x,\theta)$ satisfy the above Hypothesis \ref{hip}, and if $\delta>(\log k)^{-1}$ is such that
$$
\frac{1}{k} \frac{\phi(B)}{B} \sum_{L\in {\mathcal L}} \frac{\phi(q_i)}{a_i} \#{\mathcal P}_{L,{\mathcal A}}\ge \delta \frac{\#{\mathcal A}(x)}{\log x},
$$
then
$$
\#\left\{n\in {\mathcal A}(x): \#(\{L_1(n),\ldots, L_k(n)\}\cap {\mathcal P})\ge C^{-1} \delta \log k\right\} \gg \frac{\#{\mathcal A}(x)}{(\log x)^k \exp(Ck)}.
$$
\end{theorem}
To get to our problem, we take a large $x$ and ${\mathcal A}(x):=[x,2x)\cap {\mathbb N}$. We take $k:=\lfloor (\log x)^{1/5}\rfloor$ and
$$
L_i(n):=q^{k+1}n+q^i+1,\qquad {\text{\rm for}}\quad i=1,\ldots,k.
$$
Thus, $a_i:=q^{k+1}$ and $b_i:=q^i+1$ for $i=1,\ldots,k$. Numbers of the form $L_i(n)$ have their last $k$ digits in base $q$ equal to $0$ except for the digits in the first and $i$th position (counted from right to left) which are $1$.
Clearly, $L_i(n)\in B(m,1)$ for al $i=1,\ldots,k$, where $m:=q^{k+1}n+1$. Note that the condition $0\le a_i,b_i\le x^{\alpha}$ holds for large $x$ with $\alpha=1/3$. Hypothesis (i) and (iii) are trivially satisfied.
For hypothesis (ii), we follow the proof of Theorem 3.2 in \cite{May}.  By the Landau--Page theorem there is at most one modulus $m_0\le \exp({\sqrt{\log x}})$ such that there exists a primitive character $\chi$ modulo $m_0$
for which $L(s,\chi)$ has a real zero larger than $1-c_2(\log x)^{-1/2}$ for suitable constants $c_1,~c_2$. If this exceptional modulus $m_0$ exists, we take $B=qP(m_0)$, where $P(m_0)$ is the largest prime factor of $m_0$ and otherwise take $B=q$. If $m_0$ exists it must be squarefree apart from a possible factor of $4$ and satisfies $m_0\gg \log x$. Thus, $q\log\log x\ll B\ll q\exp(c_1{\sqrt{\log x}})$. Thus, whether $m_0$ exists or not we still have
$$
\frac{\phi(B)}{B}=\frac{q}{\phi(q)}+O\left(\frac{1}{\log\log x}\right).
$$
Now (ii) follows from a variant of the Bombieri-Vinogradov theorem avoiding an exceptional character.   Now  for each $i\in \{1,\ldots,k\}$,
\begin{eqnarray*}
\#{\mathcal P}_{L_i,{\mathcal A}}(x) & \ge & \frac{\pi(2q^{k+1}x+q^i+1;q^{k+1},q^i+1)-\pi(q^{k+1}x+q^i+1,q^{k+1};q^i+1)}{2\phi(q^{k+1})}\\
& > & \frac{q\#{\mathcal A}(x)}{3\phi(q)(\log x)}
\end{eqnarray*}
for $x$ sufficiently large (see Theorem 3.5 in \cite{BaSh}). Thus,
$$
\frac{1}{k}\frac{\phi(B)}{B}\sum_{i=1}^k \frac{\phi(a_i)}{a_i} \#{\mathcal P}_{L_i,{\mathcal A}}(x)>\frac{\phi(q)\#{\mathcal A}(x)}{4q\log x}.
$$
Hence, we may take $\delta:=\phi(q)/(4q)$. Thus, provided $k>e^{4q/\phi(q)}$, we have, by Maynard's theorem, that
$$
\#\{n\in {\mathcal A}(x): \#(\{L_1(n),\ldots,L_k(n)\}\cap {\mathcal P})>C^{-1}\delta \log k\}\gg \frac{\#{\mathcal A}(x)}{(\log x)^k \exp(Ck)}.
$$
The proof is now complete since $C^{-1} \delta \log k\gg \log\log x\gg \log \log n$ for any $n\in [x,2x)$, while $n<m=q^{k+1}n+q^i+1<n^2$ for large $x$ (so $\log\log n=\log\log m+O(1)$).

\section{Conclusion}
In this short note, we have considered the proximity of long primes under the angle of their expansion on some numeration basis. It would be of interest to generalize our coding theoretic approach to other sparse sequences of integers satisfying a prime number theorem. Number fields might provide other examples of primes enjoying similar results.
While we restricted our study to two distances, due to their favorable arithmetic interpretation, it is possible that other distances yield interesting results.
Regarding the problem of how close primes are to large integers, we had to use more classical techniques. It seems difficult to determine the covering radius of the prime code for instance,
given its large size.

\end{document}